\documentclass[letterpaper, conference]{IEEEtran}
\usepackage{graphicx, epstopdf}
\usepackage{pgfplots}
\usepackage[cmex10]{amsmath}
\usepackage{listings}
\usepackage{array}
\usepackage[caption=false,font=footnotesize]{subfig}
\usepackage{pdflscape}

\begin{document}

\title{Local Short Term Electricity Load Forecasting: \\
Automatic Approaches}

\author{\IEEEauthorblockN{The-Hien Dang-Ha\IEEEauthorrefmark{1},
		Filippo Maria Bianchi\IEEEauthorrefmark{2},
		Roland Olsson\IEEEauthorrefmark{3}}
		
		\IEEEauthorblockA{\IEEEauthorrefmark{1} Department of Informatics,
		University of Oslo, Norway, Email: hthdang@student.matnat.uio.no}
		\IEEEauthorblockA{\IEEEauthorrefmark{2} Machine Learning Group,
		University of Troms\o{}, Norway, Email: filippo.m.bianchi@uit.no}
		\IEEEauthorblockA{\IEEEauthorrefmark{3} Faculty of Computer Sciences
		\O{}stfold University College, \O{}stfold, Norway, Email: roland.olsson@hiof.no}
}

\maketitle

\begin{abstract}
\textit{Short-Term Load Forecasting} (STLF) is a fundamental component in the efficient management of power systems, which has been studied intensively over the past 50 years. The emerging development of smart grid technologies is posing new challenges as well as opportunities to STLF. 
Load data, collected at higher geographical granularity and frequency through thousands of smart meters, allows us to build a more accurate local load forecasting model, which is essential for local optimization of power load through demand side management. 
With this paper, we show how several existing approaches for STLF are not applicable on local load forecasting, either because of long training time, unstable optimization process, or sensitivity to hyper-parameters. 
Accordingly, we select five models suitable for local STFL, which can be trained on different time-series with limited intervention from the user.
The experiment, which consists of 40 time-series collected at different locations and aggregation levels, revealed that yearly pattern and temperature information are only useful for high aggregation level STLF. 
On local STLF task, the modified version of double seasonal Holt-Winter proposed in this paper performs relatively well with only 3 months of training data, compared to more complex methods.
\end{abstract}

%\begin{IEEEkeywords}
%	Short-term Load Forecasting, Demand Response, Semi-Parametric Addictive Model, Random Forest, NARX, time-series, ARMA, SARMA, DSWH, TBATS
%\end{IEEEkeywords}

%%%%%%%%%%%%%%%%%%%%%%%%%%%%%%%%%%%%%%%%%%%%%%%%
%%%%%%%%%%%%%%%%% INTRODUCTION %%%%%%%%%%%%%%%%% 
%%%%%%%%%%%%%%%%%%%%%%%%%%%%%%%%%%%%%%%%%%%%%%%%

\section{Introduction}
\label{sec:intro}

Load forecasting is an integral part of electric power system operations, such as generation, transmission, distribution, and retail of electricity \cite{Hong2015}. 
According to different forecast horizons and resolutions, load forecast problems can be grouped into 4 classes: long-term, mid-term, short-term and very short-term. 
In this paper, we focus on \textit{Short-Term Load Forecast} (STLF) of hourly electricity load for one day ahead.

Due to its fundamental role, STLF has been studied intensively over the past 50 years.
%most of the proposed methodologies are still at a theoretical level and have limited practical value \cite{Hong2014}. 
However, the deployment of smart grid technologies brings new opportunities as well as challenges to the field. 
On a smart grid, load data can be collected at a much higher geographical granularity and frequency than before, by means of thousands of smart meters \cite{Hong2015}. 
Such a larger availability of data allows the synthesis of a more \emph{local} load forecasting model, which is essential to optimize power load \emph{locally} in a demand-response paradigm \cite{chiu2009framework}. 
We refer to a load time-series as ``local'' when it contains measurements relative to a small geographical region, whose average hourly load goes from several hundreds up to several hundreds of thousands of kWh. 
 
Despite the great variety of STLF methods proposed in the literature, most of them focus on load time-series relative to high aggregation levels (big towns, cities or entire countries), whose average goes from several to hundreds of MWh \cite{Goude2014}. We found that these methods are not applicable to local STLF task for the following reasons:

\begin{itemize}
    \item \emph{Long training time:} unlike for STLF focusing on high aggregation level, in local STLF we need to train and update thousands of models at the same time. Predictions are made on hourly basis for many local regions, and the forecasting models must often be retrained (e.g., each month). This requirement cuts out approaches relying on slow derivative-free optimizers, such as evolutionary algorithms or particle swarm optimization \cite{wang2008new, wu2009novel, liao2006application}.
    
    \item \emph{Unstable optimization process: } since thousands of models needs to be trained at the same time, a local STLF model needs to be robust to the discrepancy in time-series characteristics. For example, including long-term seasonal dependency into state space models makes their optimization process unstable. 
    
    \item \emph{Sensitivity to hyperparameters: } nonparametric techniques such as artificial neural networks and kernel estimation are characterized by a high sensitivity to the hyperparameters of the model. For example, Feed-Forward Neural Networks (FNN) have been proposed and extensively used for STLF since the 1990s \cite{Alfares2002, Hong2010}. 
    However, their prediction performance highly depends on the number of layers, the amount of nodes per layer, the regularization coefficients, and the learning rate.
    Such hyperparameters must be tuned through cross-validation, which is time-consuming, due to the slow gradient descent training procedure, and does not guarantee convergence. 
    Additionally, FNN approach requires a carefully-designed preprocessing such as outliers removal to work effectively \cite{da2001enhancing, hernandez2013short}. These problems make FNN unsuitable for the local STLF task. 
    On the other hand, recurrent neural networks such as echo state networks and long short term memory networks are widely adopted in STLF \cite{7286732}.
    However, these architectures are not considered in this work, as we focus on model-based approaches.
    
\end{itemize}

After having conducted a comprehensive survey of different STLF approaches, we selected five models, two of which are original variations of existing architectures, proposed in this paper for the first time. 
The models were chosen (or modified) to overcome the three aforementioned limitations and to be at the same time characterized by a high degree of automation, both in the training and in the prediction phase. 

In our experiments, we process $40$ time-series collected from separate locations and characterized by different aggregation levels. 
To the best of our knowledge, this is the first time a local STLF experiment has been done with such a large number of time-series with different characteristics. 
As expected, our experiments show that yearly patterns and temperature information are only useful for high aggregation level STLF. 
On very local load time-series (less than several hundreds of kWh), the modified version of double seasonal Holt-Winter (modifiedDSHW) proposed in this paper performs relatively well with only $3$ months of training data, compared to other more complex methods that require years of training data.

The remainder of the paper is organized as follows. 
In Section \ref{sec:DataDescription}, the datasets under consideration are described, and the main characteristics of the load time-series are analyzed. 
The five proposed models and their origins are presented in Section \ref{sec:Methodology}. 
Section \ref{sec:Experiments} explains the experiment setup and discusses the obtained results. 
Section \ref{sec:Conclusion} concludes the paper and suggests some future work.

%%%%%%%%%%%%%%%%%%%%%%%%%%%%%%%%%%%%%%%%%%%%%%%%%%%%
%%%%%%%%%%%%%%%%% DATA DESCRIPTION %%%%%%%%%%%%%%%%% 
%%%%%%%%%%%%%%%%%%%%%%%%%%%%%%%%%%%%%%%%%%%%%%%%%%%%

\section{Data Description}
\label{sec:DataDescription}

%299 substations, >8000 smart meters
The dataset under analysis consists of $40$ load time-series collected from two countries, US and Norway, at different levels of aggregation. Such diversity in the data allows us to benchmark the generalization capability of various forecast methods. 
Among these $40$ time-series, $20$ of them come from the Global Energy Forecasting Competition 2012 (GEFCom2012). This dataset consists of 4 years of hourly load collected from a US utility with $20$ zonal level series with average hourly load varies from $10.000$kWh up to $200.000$kWh. The dataset is also accompanied by $11$ temperature time-series collected at the area, which can be used to improve the forecasting performance. 

The other $20$ load time-series come from Hvaler, a small island in Norway with around $6000$ households. The island has been used as a smart grid pilot for many years, with over $8000$ smart meters have been installed since $2012$. The island power grid includes around $100$ small distribution substations (including transformer on pole) organized hierarchically. Since there is currently no smart meters installed at these small substations, their loads are estimated by aggregating from their corresponding smart meters installed at households, street lights or other end consumers. The $20$ time-series are relative to small distribution substations. They cover two years (2012-2013) and were selected based on the quality of data (e.g. the number of missing entries). 
As opposed to the GEFCom2012 dataset, the Hvaler's load time-series are much more local and have less than $200$kWh average hourly load. This allows us to test how different predictive models perform at various aggregation levels.
Before delving into the details of each model, we first examine some characteristics of the load signals. %Performing this kind of analysis is common practice in data analytic studies for load forecasting.

\subsection{Seasonal Patterns}
Fig. \ref{fig:Hvaler_profile} shows the hourly load at Hvaler's substation 1 over 2 years.
Through a simple inspection, we can observe a strong seasonal pattern characterized by high demand for electricity in winter and low demand in summer. 
This pattern evidences a dependency between weather conditions and power consumption. 
However, such a relationship depends also on geographic location and type of consumers. 
Fig. \ref{fig:GEFCOM_profile} shows hourly load at zone 1 of the GEFCom2012 dataset, indicating high demand during summer and winter, while low demand in other seasons.
By analyzing the data more in detail, we can notice intraweek seasonal cycles (the load demand on the weekend is usually lower than on the weekdays) and intraday seasonal cycles, which arises from human routines (e.g. peaks at breakfast time and before dinner). Although these yearly, intraweek, and intraday seasonality effects are common in load time-series, their importance is not studied in local STLF. In our experiment, we observed that the yearly pattern is useful only if the load time-series is highly aggregated.

\subsection{Weather Effects}
\begin{figure*}
	\centering
	\subfloat[Hvaler - Norway\label{fig:HvalerTempCons}]{
		\includegraphics[width=.47\linewidth, trim={0cm 1.5cm 0cm 2.8cm},clip]{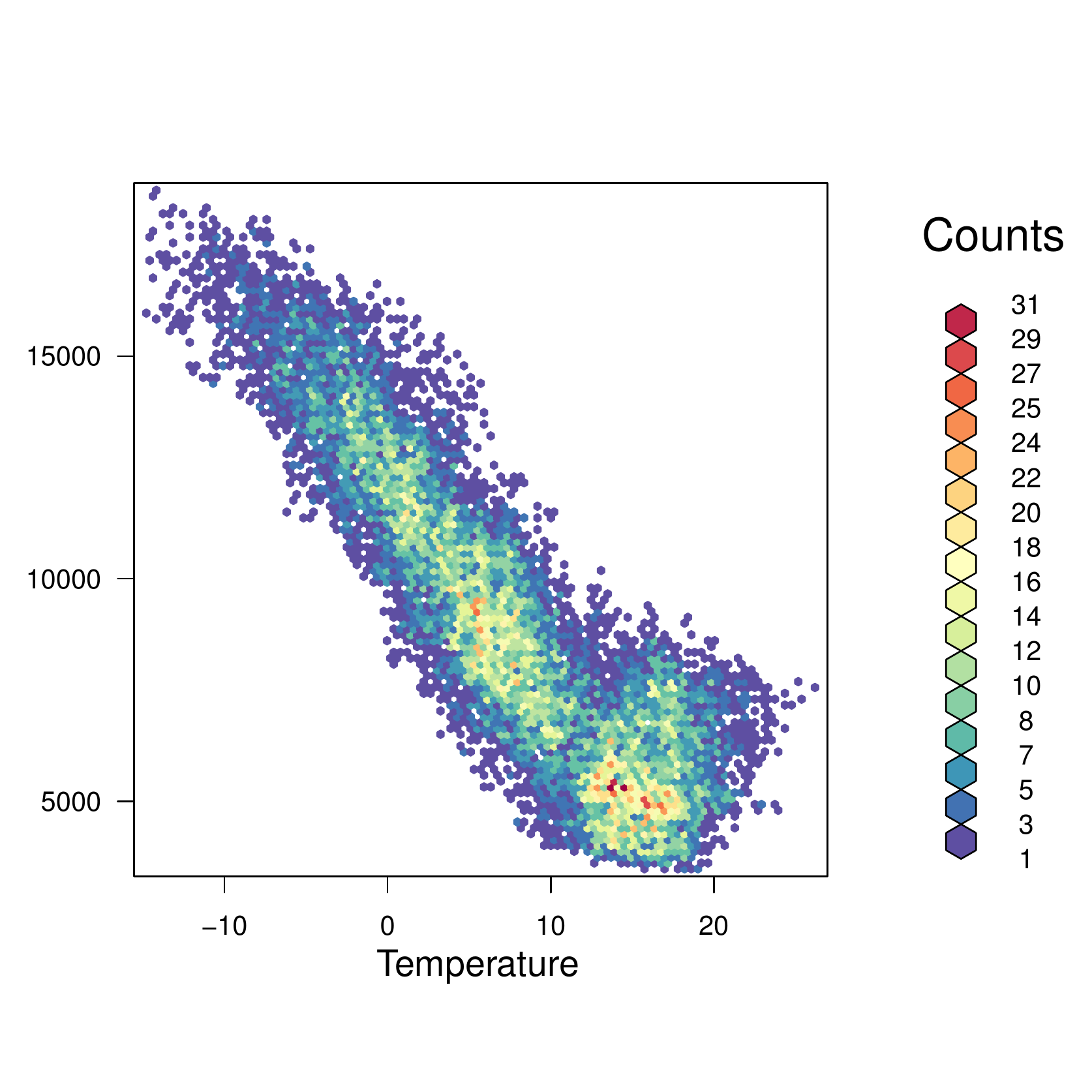}
		}
    \subfloat[GEFCom2012 - US\label{fig:Z02TempCons}]{
		\includegraphics[width=.47\linewidth, trim={0cm 1.5cm 0cm 2.8cm},clip]{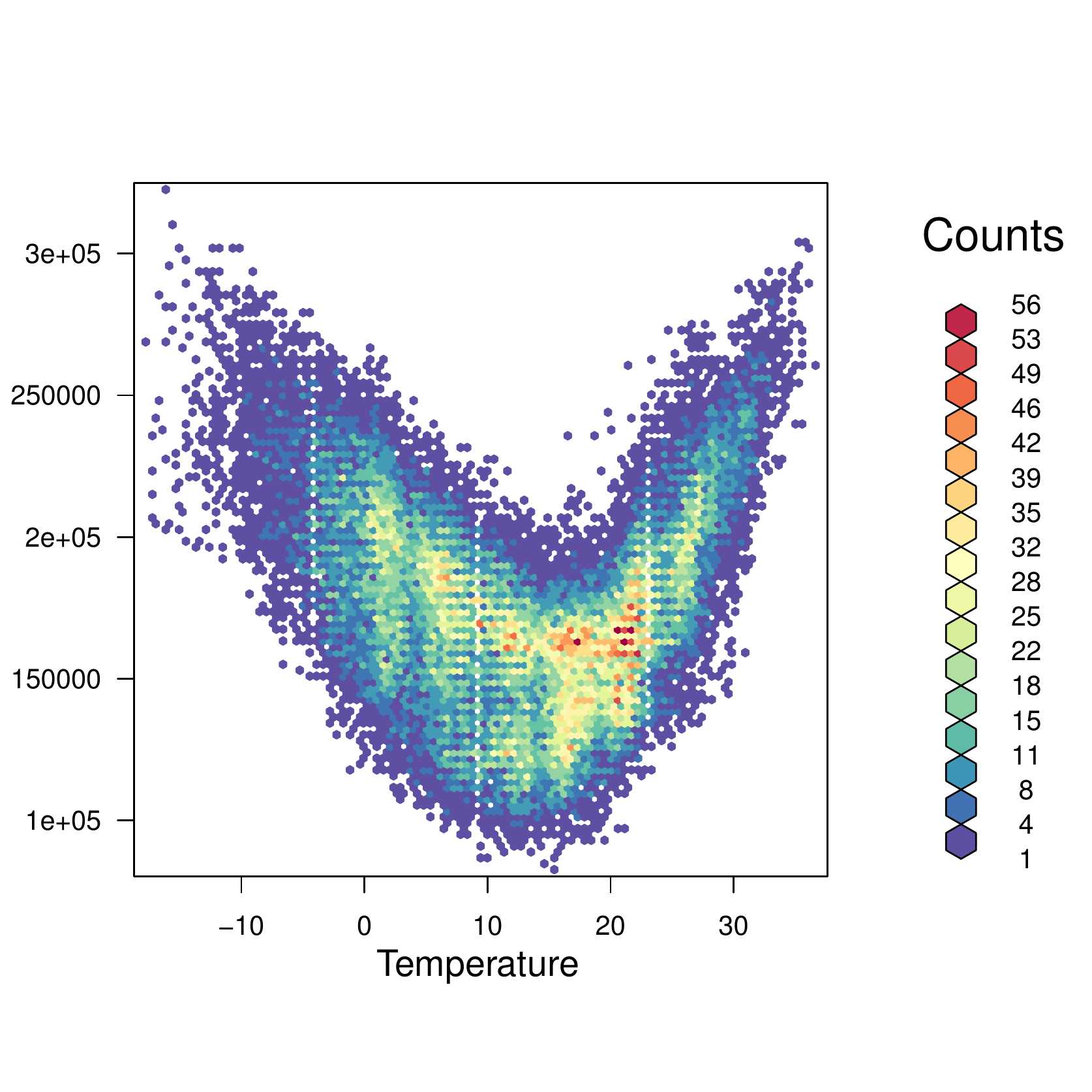}
	}
	\caption{Relationship between Temperature ($^{\circ}\mathrm{C}$) and Consumption (KWh) in Hvaler (top) and GEFCom2012 dataset (bottom)}
	\label{fig:TempCons}
\end{figure*}

In load forecasting, weather conditions have always been an important factor. 
Although many meteorological elements like humidity, wind, rainfall, cloud cover, thunderstorm could be accounted for, the most influential and popular is the \emph{temperature}, whose measurement is also easier to retrieve. 
In fact, temperature variables can explain more than 70\% of the load variance in the GEF2012Com dataset \cite{Hong2015}
The scatter plot in Fig. \ref{fig:TempCons} shows the relationship between load and temperature in both Hvaler and US. 
While the "V" shape is consistent with the two cases, there are still obvious differences on the relationship, which could be explained by the difference in geographical locations, human comfortable temperature, heating/cooling technology, or type of consumers (e.g. industrial or residential units).

\subsection{Calendar Effects}
People change their daily routines on calendar events, such as holidays, festivities and special events (e.g. football matches, transportation strike), with a possible modification in the electricity demand. 
Those situations represent outliers and could be treated differently to improve the model accuracy.
In our study, only the national holiday events are taken into account, due to the lack of information on other events. 

\subsection{Long-Term Trend}
\begin{figure*}
\centering
    \subfloat[Load profile of substation 1 from Hvaler with the 4 testing periods are masked out.]{
    \includegraphics[width=0.7\textwidth,trim={0.7cm 4cm 1.2cm 4.3cm},clip]{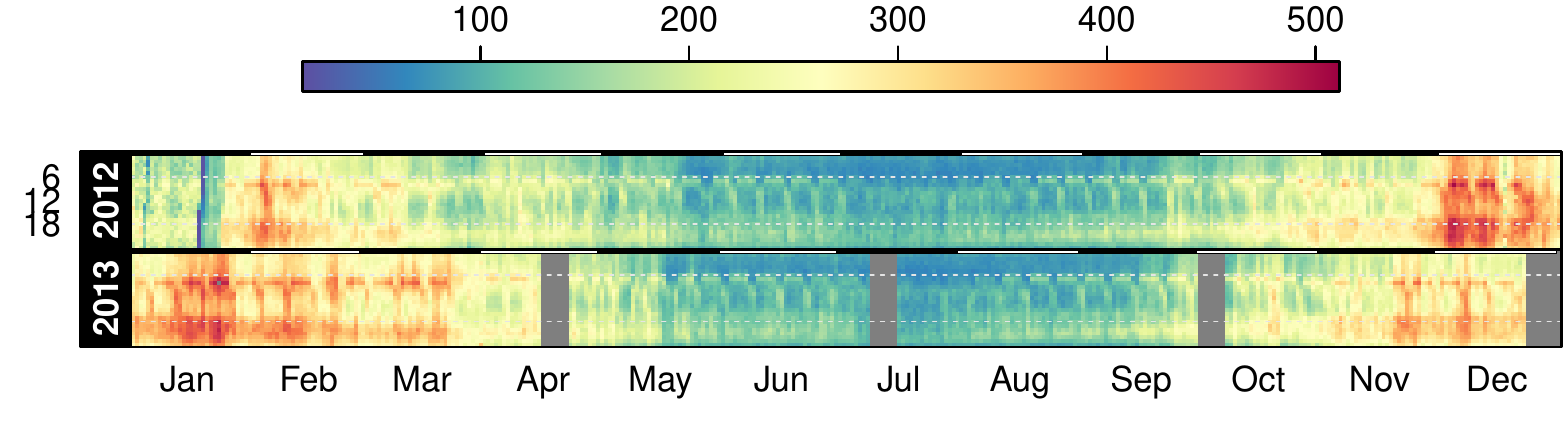}
    \label{fig:Hvaler_profile}}
    
    \subfloat[Load profile of zone 1 from GEFCOM2012 (in kWh) with the 4 testing periods are masked out.]{
    \includegraphics[width=0.7\textwidth,trim={0.7cm 3cm 1.2cm 3.3cm},clip]{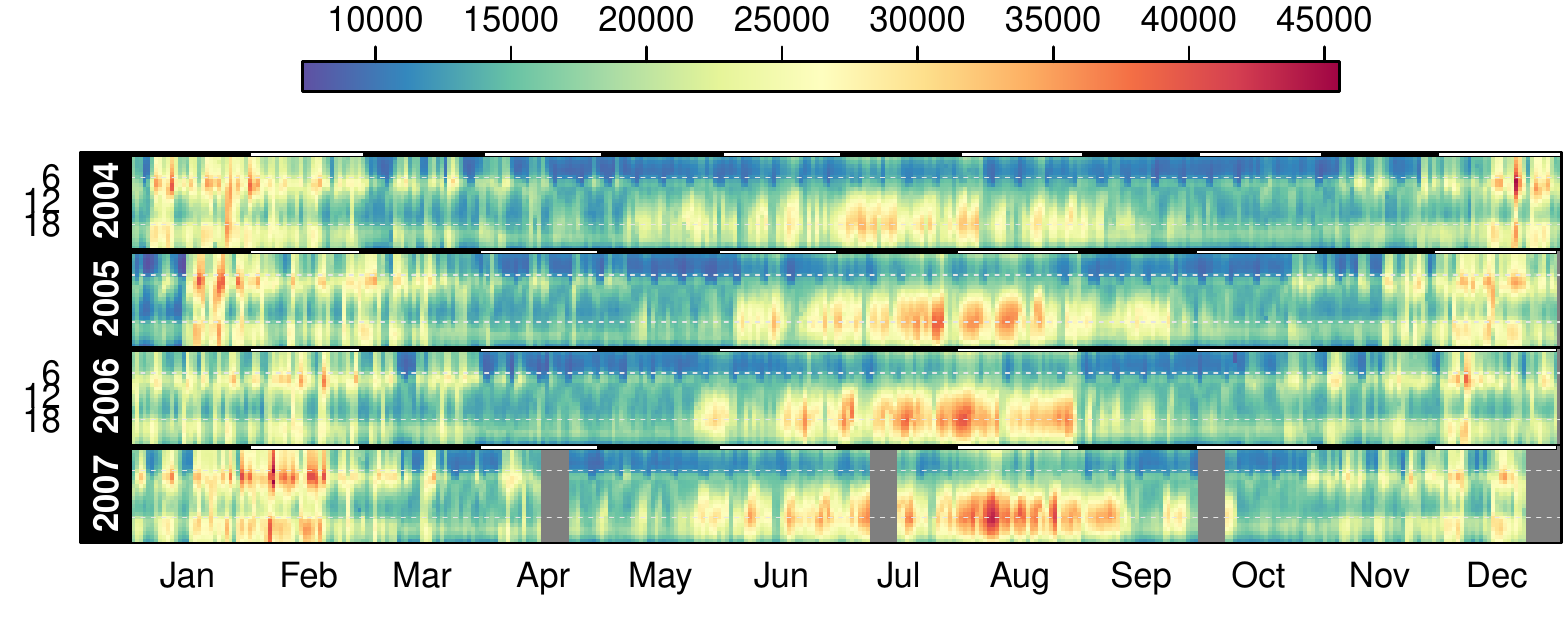}
    \label{fig:GEFCOM_profile}}
\caption{Hourly load profiles in KWh of the substation 1 from Hvaler (top) and zone 1 from GEFCom2012 dataset (bottom).}
\label{fig:Consumption}
\end{figure*}

The scale, variation, and other properties of the load signal could change over time due to changes in population, technology, or economic conditions. In Fig. \ref{fig:Consumption}, we can see a tendency of increasing consumption over the years. In our experiment, two methods explicitly model and detrend the load time-series in the first step. Other methods either model it implicitly or ignore it. This is because the long-term trend can be considered as constant in short-term and may not contribute much to one-day-ahead load forecasting.

%\subsection{Load Decomposition Approach}
%As a summary of the above observations, load demand could be considered as a signal arising from multiple underlying processes:
%\begin{itemize}
%	\item Daily, weekly, and yearly routines: people begin and end the day and vary their energy consumption in a regular manner. However, these routines could change over time.
%	\item Special people activities on special events.
%	\item Non-human consumption, such as street lighting or industrial consumers.
%	\item Consumption affected by weather deviation from seasonal norms.
%	\item System scale changes over time due to changes in population, economic conditions, or technology
%\end{itemize}

%Decomposition provides a basis for us to analyze the abilities of different forecasting models on capturing the effects of different processes, as well as diagnose their limitations for further improvement.

%\commentF{In this last subsection is not very clear what you want to say. If the content of this last subsection is not necessary to understand the rest of the paper, maybe it could be removed? }


%%%%%%%%%%%%%%%%%%%%%%%%%%%%%%%%%%%%%%%%%%%%%%%
%%%%%%%%%%%%%%%%% METHODOLOGY %%%%%%%%%%%%%%%%% 
%%%%%%%%%%%%%%%%%%%%%%%%%%%%%%%%%%%%%%%%%%%%%%%

\section{Methodology}
\label{sec:Methodology}
In this section, we review several important STLF approaches. We divided the STLF methodologies of interest into four main categories: \emph{Averaging}, \emph{Linear State Space}, \emph{Decomposition}, and \emph{Data-Driven} approaches. 
This classification is particularly meaningful for our analysis, but other taxonomies have been proposed in the literature \cite{Hong2015, Alfares2002, Aman}. 
For each category we provide a short overview, pointing out the main advantages and disadvantages encountered in our applicative scenario. 
Based on this analysis, one baseline and five potential solutions were implemented and tested on the Hvaler and GEFCom2012 data. 

% ------------------ Averaging Approach ------------------

\subsection{\textbf{Averaging Approach}}
Despite being very basic, this is still a popularly used method due to its simple implementation and obvious model interpretation \cite{coughlin2008estimating}. 
The averaging model makes predictions based on linear combinations of consumption values from ``similar'' days and was used as a benchmark in \cite{Taylor2007}. The forecast is computed as
\begin{equation}
\label{eq:average}
\hat{y}_t(k) = \frac{y_{t + k-s_2} + y_{t + k-2s_2}+y_{t + k-3s_2}+y_{t + k-4s_2}}{4},
\end{equation}
where $y_t$ is the demand in period $t$, $k$ is the forecast lead time, and $s_2$ is the so called second seasonal cycle, which is the intraweek $s_2 = 24*7 = 168$. 
The model predicts future load by averaging the corresponding observations in each of the previous four weeks. As previously discussed, there are 3 typical seasonal cycles in load time-series, which are intraday, intraweek and yearly cycles. In this paper, their cycle lengths (in hour) are represented respectively as $s_1 = 24$, $s_2 = 24*7=168$, and $s_3=24*364.25=8766$. Note that we fixed the length of each seasonal cycle and did not re-estimate it in different time-series. This is because the estimation requires analyzing a periodogram manually and is hard to automate.

At first, we intended to use the averaging model as the baseline to estimate the forecast difficulty. However, the model produced highly autocorrelated errors in the residual series and could not serve as a good baseline. 
Therefore, we decided to train an additional ARIMA model on the residual series produced by the averaging model.
We call this the \emph{avgARIMA} model and use it as the baseline model in our experiment.

% ------------------ Linear State Space Approach ------------------

\subsection{\textbf{Linear State Space Approach}}
State space approach refers to those models that can be written in a linear state space form, which consists of a set of states with initial distribution (usually Gaussian), a measurement equation and a Markovian transition equation. 
Although state space models can be extended to include exogenous variables, such as temperature, the univariate setting is still the most popular in STLF. 
Recent studies have shown that, although in the long run the load is strongly influenced by meteorological conditions and special events, an univariate model is sufficient in shorter lead times \cite{Taylor2006}. 
%In fact, Taylor et.al. \cite{Taylor2006, Taylor2007, Taylor2012} have shown that univariate models can yield superior performance for lead times up to six hours ahead. 
The two most common and accurate state space models that have been reported for STLF task in the literature are \textit{Auto-Regressive Integrated Moving-Average} (ARIMA) and \textit{Holt-Winters exponential smoothing}. 

\subsubsection{ARIMA}
The ARIMA model was adopted in STLF back in 1987 \cite{Hagan1987}, where a double seasonal ARIMA model was experimented (intraday and intraweek cycles). 
This approach remains popular, with extensions to include exogenous variables, or intrayear seasonal cycles \cite{wang2008new, Taylor2010}. 
One big disadvantage of ARIMA is that the model hyperparameters (such as AR, I, and MA orders, as well as orders of seasonal AR, I, and MA terms) are usually derived from the Box-Jenkins test, which is hard to automate and still require human expertise to examine the partial correlogram of the time-series \cite{Ngo2013, Aman}. 
These hyperparameters can be heuristically fixed \textit{a-priori} \cite{Taylor2007, Taylor2010}, However, during our experiment we noted that in the ARIMA model with (double) seasonality the optimization process becomes highly unstable when fixing the AR, I, MA orders and all the seasonal AR, I and MA orders to some arbitrary numbers. %\commentF{Move this in the experiments/discussion -- We consistently got convergence error on some random set of time-series even though we have tried various combinations of different optimizer and state initialization scheme. This happened to the double seasonal ARMA model \cite{Hagan1987} and the triple seasonal Holt-Winters model \cite{Taylor2010}, which are arguably the two most accurate state space models that has been reported for STLF task in the literature.}
%The optimizer failed to converge in some time-series after a very long running time, even though various optimizers and initialization schemes have been tried.
Akaike Information Criterion can also be used to set ARIMA hyperparameters. However, a complete search over all possible models is time-consuming, especially for seasonal ARIMA. 
Therefore, in this paper, we only use ARIMA to correct the autocorrelation in the residual series produced by other models. This is a common practice in time-series forecasting to improve accuracy when the main model has autocorrelated errors.

\subsubsection{Holt-Winters}
Holt-Winters is another popular state space model that accommodates the two intraday and intraweek seasonal cycles that commonly appear in load time-series. 
Taylor et.al. (2003) \cite{Taylor2003} introduced the double seasonal Holt-Winters method (DSHW), whose important advantage, which we find suitable for our local load forecasting problem, is that it only requires the length of the two seasonal cycles to be specified. 
We indeed did not get any optimization problem when fixing these two seasonal cycles to the intraday $s_1 = 24$ hours and intraweek $s_2 = 168$ hours. 
Implementation details are provided in Section \ref{sec:Experiments}, where we also introduce a modified version of the original DSHW that yields significantly better performance in our experiment.

In 2010, Taylor et.al. \cite{Taylor2010} proposed the triple seasonal Holt-Winters (TSHW), which enables the model to accommodate the intrayear seasonal cycles. 
However, it turned out that the extra seasonality causes the training process to become much slower and unstable. In fact, the success of the optimization process depends on the choice of initial values for the initial states. 
To address this issue, Taylor et.al. generated $10^4$ initial vectors as possible initialization for the variables of the model. Since this process requires significant computational power, the TSHW is unsuitable for the local STLF task.
%The model is available through the \emph{forecast} package in $R$ under the name \emph{dshw}. TODO: Move

% ------------------ Data-Driven Approach ------------------

\subsection{\textbf{Data-Driven Approach}}
Instead of modeling the underlying physical processes, data-driven methods try to discover consistent patterns from historical data, according to a machine learning approach. 
A mapping between the input variables and the load is learned and then used for prediction. 
Depending on the forecasting task, the input and output variables are designed accordingly.
In the following, we review two of the main approaches in STLF.

\subsubsection{Nonlinear Non Auto-Regressive Regression} 
\label{Non-linear Regression}
This approach models the load as a non-linear function of only exogenous input variables without using the autoregressive terms.
According to the nature of the data discussed in Section \ref{sec:DataDescription}, potential exogenous variables of interest are: time of day, time within week, time of year, linear trend, temperatures, or smoothed temperatures. 
Frameworks such as random forests and gradient boosted model\dots can be used to map inputs to the desired output. 
%In GEFCom2012, the second winning team was using gradient boosted model to build a non-linear regression model as part of their solution, using time of day, time within week, temperatures, and smoothed temperatures as the input variables \cite{Lloyd2014}. 
However, since these approaches do not explicitly model the autocorrelation that almost always exists in the load time-series, they must be used in combination with other techniques to be effective, such as a state space and a long-term decomposition model \cite{Lloyd2014}. 
This requires extra effort on the model deployment and management process. For this reason, we did not include this approach in our study.

\subsubsection{Nonlinear Auto-Regressive with Exogenous (NARX)} 
%This model use non-linear learner, such as ANN or SVM, to learn the patterns between future load and its previous values together with exogenous variables 
An NARX model computes the next value of a variable, from previous values of the variable itself and current and past values of exogenous series. 
The basic formulation reads as
\begin{equation}
    \label{eq:narx}
    y_{t}=F(y_{{t-1}},y_{{t-2}},\ldots ,x_{{t}},x_{{t-1}},x_{{t-2}},\ldots ),
\end{equation}
where $F(\cdot)$ is a non-linear function, which could be modeled by any general-purpose machine learning model such as artificial neural network (ANN) or support vector machine (SVM). 
%This can be referred as \emph{Time-Delay Network} if an ANN is used.

%In \cite{Taylor2012}, a similar model is tested, which yields comparable performance to other improved versions of DSWH. We found this model appropriate for our problem, thanks to its small training data requirements, easy implementation, and good performance. 
As discussed in Section \ref{sec:intro}, ANN depends on several hyper-parameters and its training can be cumbersome.
An ANN is prone to overfitting and is sensitive to outliers \cite{Hong2015}. %Hence, if the load signal is noisy or extreme weather condition happens, the model will produce bad predictions.
The problem could be solved by replacing the ANN with a model characterized by a lower variance such as SVM, which has been adopted in several studies \cite{Chen2001,Aung2012,Sapankevych2009}. 

A different approach is to use the random forest (RF) as in \cite{breiman2001random}. 
Thanks to its bagging, data sub-sampling, and random features selection process, the RF model is capable of capturing complex patterns, while maintaining a low variance \cite{liaw2002classification}. 
This approach is called \emph{NARX-RF} and is specified in Section \ref{sec:specification}.

%\commentF{Subsections C.1 and C.2 are very similar. Maybe merge them into a single section called data-driven approaches? Also, non-linear regression is a really wide area, maybe it is safer to avoid to use that word :) Dang: I Can't find an easy way to merge :(, so I just change its name to non auto-regressive  }

% ------------------ Time-Series Decomposition Approach ------------------

\subsection{\textbf{Time-Series Decomposition Approach}}
\label{Decomposition}
Time-series decomposition approach deconstructs a time-series into several components, each component represents different kinds of pattern. According to the nature of load time-series discussed in Section \ref{sec:DataDescription}, potential components are long-term trend, intrayear cycle, intraweek cycle, intraday cycle, relationship between temperature and load, holiday events, \dots

In GEFCom2012, Lloyd and James (2014) \cite{Lloyd2014} used a \emph{Gaussian Process} to decompose the load time-series. Their Gaussian Process contains a set of different kernels, each of them is designed to captured different component in the time-series. The long-term trend is captured by a squared exponential kernel, intrayear cycle is captured by a periodic kernel of time, while the relationship between temperature and load is captured by a squared exponential kernel of temperature. Although this hybrid approach works relatively well on GEFCom2012 data, a GP model needs to be manually carefully designed and requires special treatment on different load signals \cite{Lloyd2014}. Therefore, we found Gaussian Process not suitable for the local STLF task and did not include it in the experiment.% Furthermore, GP model needs at least two years of training data. Despite these facts, we did try experimenting with a GBM and GP combination model on the GEFCom2012 dataset, since the implementation is generously made available by the authors of \cite{Lloyd2014}. However, the \textit{Mean Absolute Percentage Error} (MAPE)\footnote{MAPE is a widely used benchmark in energy forecasting community}, rolling testing on the last month June-2018 at Zone 2, is $9.5\%$, which is poor \commentF{Also this detail should go in the experimental section}
%\footnote{The original model also includes another linear regression model, which we find from their source code that it does not contribute much to the result} compared to our other proposed methods. Of course, we could  improve the model by using a time-series model like ARMA to capture auto-correlation in the residual time-series. 
%However, we found that with more than two years of training data, a semi-parametric addictive model is more appropriate for the task and can bring better result (see section \ref{LinearRegression}).

%-------Linear Regresion---------
Another popular way to decompose a load time-series is to use linear additive models, where the load is modeled as a linear combination of various independent features. The learned model is interpretable, easy to implement/automate the training process, and able to achieve high accuracy. %This statement is partly proven at the GEFCom2012 competition, where all of the three best solutions has a linear additive model \cite{Hong2014}. 
Many different features have been suggested in the literature to capture different load components. For example, the yearly cycle can be modeled by $8$ Fourier series \cite{Dordonnat2008} or spline functions \cite{Nedellec2014}; relationship between temperature and consumption can be modeled by a pice-wise linear \cite{Dordonnat2008}, quadratic \cite{Charlton2014}, or spline functions \cite{Nedellec2014}; monthly change in relationship between temperature and consumption can be modelled by interaction terms between temperature and month of year variable \cite{Hong2010Interaction}. Among different linear additive models proposed in the literature, we found that the TBATS model suggested by De Livera et. al. 2011 \cite{de2011forecasting} and the semi-parametric additive model suggested by Goude et.al. 2014 \cite{Nedellec2014} are the two most suitable models for the local STLF task.
%However, model performance depends heavily on the choice of explanatory variables.
%We believe with enough data and well-designed features, multiple linear regression would yield the best accuracy. This statement is partly proven at GEFCom 2012 competition, where all the three best solutions include linear regression model \cite{Hong2014}.

The name TBATS is an acronym for key features of the model: Box-Cox transforms, ARMA errors, Trend, and Trigonometric Seasonal components. On the other hand, the semi-parametric additive approach bases on the Generalized Addictive Model. Precise specification of these two models are given in Section \ref{sec:specification}.

\section{Experiments and Discussion}
\label{sec:Experiments}
\label{MAPE}
\subsection{Experiment Setting}
%About the forecasting procedure, we have 4 testing periods in total, which are the four following weeks of the final year: 16, 28, 40, 52, 53 (Week 53 is just 1 day or no day at all). These are chosen to benchmark models better for different seasons. For each testing period, I do “Muti-step rolling forecasts without re-estimation”. For different testing period, I re-estimate the model.
 
%You can check this page for more detail: http://robjhyndman.com/hyndsight/rolling-forecasts/
 
In this section, we show experiment results of the five chosen models on $40$ time-series described in section \ref{sec:DataDescription}. Each time-series is marked with $4$ testing periods, which are the four following weeks of the final year data: $16$, $28$, $40$, $52$, $53$ (note that week $53$ contains just $1$ day or no day at all). These weeks were chosen to bring a fair estimation of model performance in different seasons and holidays during the final week of the year. These testing periods are demonstrated in Fig. \ref{fig:Consumption}. For each testing period, we did \emph{multi-step rolling forecasts without re-estimation}. For different testing period, we re-estimated the model. We did not re-estimate the models within one testing period since this approach is impractical. Retraining and updating thousands of models every hour or even every day are infeasible since they require a tremendous amount of computing resource. We used the \emph{Mean Average Percentage Error} (MAPE) to compare model performance, which is widely used in the energy forecasting community. The MAPEs are calculated separately at 24 prediction horizons for 24 hours ahead. Besides accuracy, we also report the training time of each model, which is an important factor in deciding which model to use in practice but rarely mentioned in the literature. Precise model specifications of one baseline and all five chosen methods are given in the following section. The whole experiment can be easily reproduced from the data and code publicly available at: \texttt{https://github.com/Nikasa1889/R\_Notebooks}

%If we re-estimate the model multiple times within one testing period, the performance can be improved significantly for some models.

%An overview of the five chosen models are presented in Table \ref{modelOverview}, 
\subsection{Model Specifications}
\label{sec:specification}
\subsubsection{avgARIMA}
\label{sec:avgARIMA}
The \emph{avgARIMA} model was used to give a baseline performance for our experiment. First, an averaging model is used to predict the future load by averaging the corresponding observations in each of the previous four weeks, as specified in (\ref{eq:average}). Its $3$-month residual time series is then used to train an ARIMA. A stepwise search was used to optimize the AR, I and MA orders. The procedure is implemented by the \emph{auto.arima()} function in the R \emph{forecast} package.

\subsubsection{originalDSHW and modifiedDSHW}
The multiplicative formulation for the original double seasonal Holt-Winders model (\emph{DSHW}) is given in the following expression \cite{Taylor2003, Taylor2007}:

\begin{equation}
l_t = \alpha(y_t/(d_{t-s_1} w_{t-s_2})) + (1-\alpha)l_{t-1}
\end{equation}
\begin{equation}
d_t = \theta(y_t/(l_tw_{t-s_2})) + (1-\theta)d_{t-s_1}
\end{equation}
\begin{equation}
w_t = \omega(y_t/(l_td_{t-s_1})) + (1-\omega)w_{t-s_2}
\end{equation}
\begin{equation}
\hat{y}_t(k) = l_td_{t-s_1+k}w_{t-s_2+k}+\phi^{k}(y_t - (l_{t-1}d_{t-s_1}w_{t-s_2}))
\label{dshw_forecast}
\end{equation}

, where $l_t$ is the smoothed level; $d_t$ and $w_t$ are the seasonal indices for the intraday and intraweek seasonal cycles, respectively; $\alpha$, $\theta$, and $\omega$ are the smoothing parameters. The term involving the parameter $\phi$ in the forecast equation (\ref{dshw_forecast}) is a simple adjustment for first-order autocorrelation. This model is implemented in the R \emph{forecast} package under the name \emph{dshw()}, and named \emph{origDSHW} in this paper.

During the experiment, we recognized that performance of the origDSHW model can be improved significantly if we employ a different objective function. Instead of using the sum of square errors of the in-sample 1-hour ahead forecast, we modified the objective function to the sum of squared errors of all 24 horizons in-sample forecast. Moreover, we also increased the upper limit of the $\phi$ parameter from $0.9$ to $0.99$. This model is called \emph{modDSHW} in this experiment. Both the origDSHW and modDSHW models were trained on $3$ months of data.
\subsubsection{NARX-RF}
\label{sec:narxRF}
Although an NARX model with SVM yielded good performance in other studies, we found it hard to automate, since its performance depends heavily on the choice of hyper-parameters: the cost of errors $C$ and width of the $\epsilon$-insensitive tube. Moreover, optimal values of these hyper-parameters vary very much on different load signals. Therefore, instead of SVM, the random forest was used, which is referred to as the \emph{NARX-RF} model in this paper. For building the random forest, the package \emph{ranger} in R was used.

To avoid multi-step ahead predictions, a separate random forest was used for each lead time. For lead time $h$, the set of inputs consists of the load values at lags: $1$, $2$, $3$, $s_1-h$, $2s_1-h$, $3s_1-h$, $s_2-h$, $2s_2-h$; two temperature-related exogenous variables: temperature and exponential smoothed temperature; and three calendar variables: time of day and day of week. The smoothed temperature is often used in STLF to take into account the physical inertia of buildings and delay effects of temperature on consumption \cite{Dordonnat2008}. The coefficient for the temperature exponential smoothing process was set to $0.85$.

We kept all the default settings in the \emph{ranger} function, which set the number of trees to grow $ntree = 500$, and the number of candidate features at each split $mtry = 3$. The subsampling ratio was set so that each tree receives $5000$ data points to train on. The model makes use of all the available data up to the testing point.
\subsubsection{TBATS}
\label{sec:TBATS}
The \emph{TBATS} model was introduced by De Livera et.al. in 2011 \cite{de2011forecasting} to solve the forecasting problem in time series with complex seasonal patterns such as multiple seasonal periods or high-frequency seasonality. The model incorporates Box-Cox transformations, linear trend, Fourier representations with time-varying coefficients, and ARMA error correction. The method involves a simple, yet efficient estimation procedure, which makes it suitable for the local STLF problem. In this experiment, the exact TBATS model described in \cite{de2011forecasting} was used without any modification. The TBATS implementation provided in the \emph{forecast} package was used with all the default settings unchanged.

\subsubsection{SemiParametric}
\label{sec:SemiPar}
The semi-parametric additive model was first introduced by Goude et.al. in the GEFCom2012 competition \cite{Nedellec2014}. In 2014, Goude et.al. have tested the method's generalization ability, where it was used for short and medium-term load forecasting on $2206$ large-scale substations automatically \cite{Goude2014}. Here we present a short explanation of the method, together with some small adaptation we have done to make it more appropriate for the local STLF task. For short, this method is called \emph{SemiPar} in this paper.

The SemiPar method splits the load into three parts:
\begin{equation}
Z_t = Z_t^{lt} + Z_t^{mt} + Z_t^{st},
\end{equation}
where $Z_t$ is the electrical load at time $t$, $Z_t^{lt}$ is the long-term part of the load, corresponding to low-frequency variations such as long-term trends or economic effects. $Z_t^{mt}$ is the medium-term part, incorporating daily to weekly effects, the meteorological effects, and the calendar effects. 
The short term part, $Z_t^{st}$, contains everything that could not be captured on a large temporal scale but could be obtained locally in time. We implemented the $Z_t^{lt}$ and $Z_t^{mt}$ exactly the same as described in \cite{Nedellec2014}. However, for the short term $Z_t^{st}$ part, we use an ARIMA model (optimized by the auto.arima() function) to capture the auto-correlation in the residual time-series after removing the long-term and medium-term parts.

The long-term forecast uses combination of generalized additive models (GAM) and kernel regression, while the medium-term forecast uses GAMs. 
The GAM model with generalized cross validation criterion is implemented in the R package  \emph{mgcv}, while the kernel regression is the Nadaraya-Watson model, which is available through the \emph{bbemkr} package.

For the long-term model, we aggregate the consumption and temperature by month, denoted by $Z_t^{monthly}$ and $T_t^{monthly}$. Then we estimate monthly consumption using the following semi-parametric additive model \cite{Nedellec2014}:
\begin{equation}
\hat{Z}_t^{monthly} = \sum_{q=1}^{12} c_q I_{Month_t = q} + f(T_t^{monthly}) + \epsilon_t
\end{equation}
Where:
\begin{itemize}
	\item $I_{Month_t = q}$ is an indicator variable which is equal to 1 when the month at observation $t$ is $q$ (from 1 to 12), and 0 otherwise.
	\item $f$ is the effect of the monthly temperature, estimated by thin plate regression splines (default setting in \emph{mgcv} package).
\end{itemize}
The monthly estimated residuals are then obtained as follows:
\begin{equation}
\hat{\epsilon}^{monthly} = Z_t^{monthly} - \hat{Z}_t^{monthly}
\end{equation}

Then the residuals are smoothed and interpolated to hourly frequency by using Nadaraya-Watson kernel regressors, with Gaussian kernels and a bandwidth of $12$. These smoothed residuals are a good estimate of low-frequency effects, which contains neither annual seasonality nor weather effects. These residuals are considered as $Z_t^{lt}$ and smooth by construction, and thus they are constantly extrapolated for one-day horizon.

By removing $Z_t^{lt}$ from the original load, we get the signal $Z_t^{det}$ which contains $Z_t^{mt}$ and $Z_t^{st}$. We fit one mid-term model for each hour of the day so that we have 24 mid-term models. These mid-term models are GAM in the following form \cite{Nedellec2014}:
\begin{align*}
Z_t^{det} = \sum m_q I_{DayType_t = q} + g_1(\theta_t) + g_2(T_t) \\ + h(toy_t) + \epsilon_t
\end{align*}
where: 
\begin{itemize}
	\item $Z_t^{det}$ is the de-trended electrical demand at time $t$.
	\item $DayType_t$ is type of day for observation $t$. 1 for Sunday, 2 for Monday, 3 for Tuesday, 4 for Wednesday, 5 for Thursday, 6 for Friday, 7 for Saturday, 8 for Christmas and New Year’s Day, 9 for Christmas Eve, 10 for Independence Day, and 11 for Thanksgiving.
	\item $\theta_t$ is the smoothed temperature, obtained via exponential smoothing of the real temperature $T_t$: $\theta_t = (1-0.85)T_t + 0.85\theta_{t-1}$
	\item $toy_t$ is the time of year, which is the position of the observation $t$ within the year. $h(toy_t)$ corresponds to the smooth yearly cycle of the load.
	
	\item All the $g(.)$ functions are modeled by thin plate regression splines, while the $h(.)$ function is modeled by cyclic cubic regression splines.
\end{itemize}

An ARIMA short term model is then built to capture patterns in the residuals after removing $Z_t^{lt}$ and $Z_t^{mt}$ from $Z_t$.

%As described above, the semi-parametric addictive approach has only two hyper-parameters: the bandwidth for kernel regression, and the coefficient for temperature exponential smoothing. However, the prediction results are quite insensitive to these parameters around the optimal values, as also justified in \cite{Nedellec2014}. The model took third place in GEFCom2012, but it gives the best result if actual future temperature data is provided. And because accurate temperature forecasts for 24-hour ahead are usually widely available through third-party services, we believe that this model is the best candidate for our local short term load forecasting task. The only issue we found when applying this method is that it needs at least two years of training data. This could be overcome by using the DSHW or NARX approaches proposed above.

\subsection{Experiment Results}
\begin{figure}
\centering
\includegraphics[width=\columnwidth,trim={0cm 0.1cm 0.2cm 0.2cm},clip]{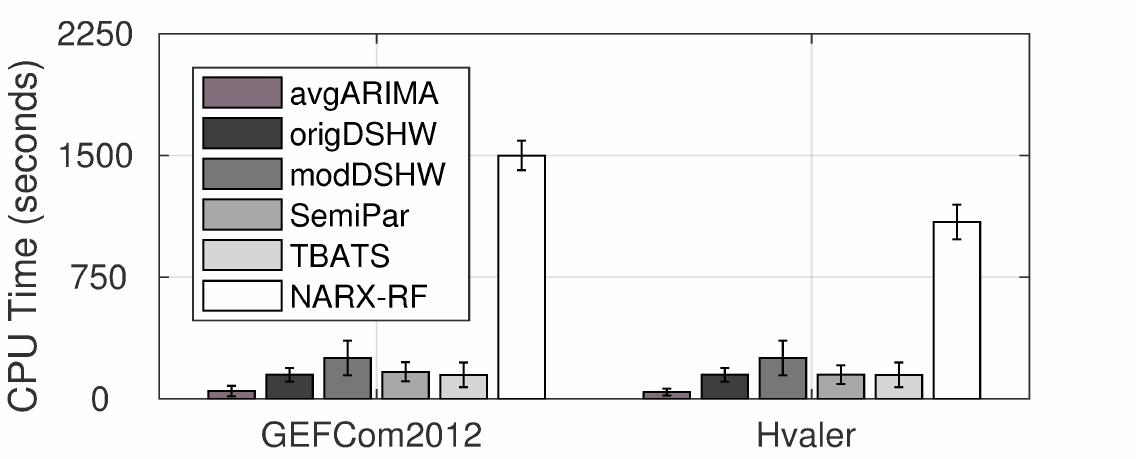}
\caption{Training time of each method running on an Intel core i7-6700k 4.0Ghz 8 cores. The CPU time was measured for one core.}
\label{fig:Time}
\end{figure}

The whole experiment was trained on an Intel Core i7-6700k 4.0Ghz machine with 8 cores. Training time of each method is reported in Fig. \ref{fig:Time}, where the CPU time was measured for one core. The whole experiment took about 12 hours to complete. 

\begin{figure}
\centering
\includegraphics[width=\columnwidth,trim={0cm 0cm 0cm 0cm},clip]{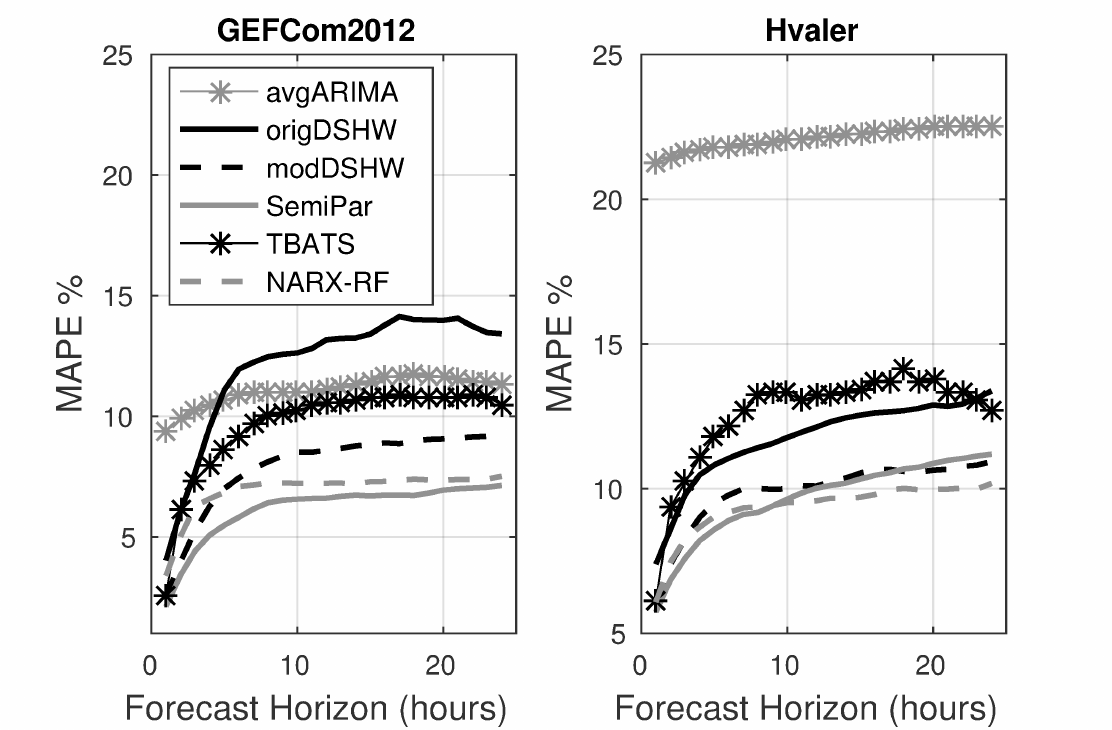}
\caption{Comparison between median MAPE of each method at different prediction horizons on GEFCOM2012 (left) and Hvaler (right)}
\label{fig:HorizonAcc}
\end{figure}

Fig. \ref{fig:HorizonAcc} shows a comparison between median MAPE of each method at different prediction horizons on the two datasets GEFCom2012 and Hvaler. On the GEFCom2012 dataset, the SemiPar method is obviously the best model. This is expected since it is the only model that has been tested and performed well in a large-scale experiment with thousands of time-series. It is also the only model that explicitly captures all the patterns discussed in section \ref{sec:DataDescription}, including the long-term, mid-term, and the short-term patterns together with the temperature effect. However, on the Hvaler dataset, where all the loads are collected at a much lower aggregation level, the SemiPar method exhibits its limitation. It performs only slightly better than the NARX-RF approach in the first ten horizons and then becomes worse when the prediction horizon increases. This can be explained by the fact that the load time-series in Hvaler are much noisier than in GEFCom2012 dataset since they consist of only a small number of consumers. Therefore, their long-term and mid-term trends are less consistent in the long run, which causes the decomposition approach to be less effective. The load in Hvaler is also collected in a shorter period (2 years), which causes the estimation of the long-term and mid-term components to become less accurate. On the other hand, the short-term processes like intra-week, intra-day, and innovations become more influential in Hvaler dataset. This explains why the modDSHW method, which only uses $3$ months of training data and only models the intra-week and intra-day seasonality, can slightly outperform the SemiPar at horizons further than $10$. 

The second best model in both of the two GEFCom2012 and Hvaler datasets is NARX-RF. It performs consistently well on the two datasets at all the prediction horizons, only significantly worse than SemiPar for early horizons. This consistency in performance is an important advantage if the system contains load time-series collected from many different scales, and we want to use only one forecasting method to simplify the deployment process. However, one has to consider its running time, since NARX-RF is more than one order of magnitude slower than other methods.

Our modifications made for the DSHW model turn out to be very effective. The modDSHW method significantly outperforms the orgiDSHW method in every case. This suggests that one should always follow these modifications if the DSHW model is of his interest. The modDSHW performs surprisingly well on the Hvaler dataset, even without using temperature information and was trained in only $3$ months of data. Therefore, we believe that temperature information only contributes marginally to the forecasting accuracy on very local load time-series like Hvaler. 

The TBATS method performs badly on both datasets and even worse than the averageARIMA model at some points. This can be because the way it decomposes the time-series is not suitable for the load signal.

%From these observations, we suggest that for 24 hours ahead forecasting, a SemiParametric-ARMA model should be %considered first (when more than 2 years of training data and temperature data are available). In case of limited historical data, or if we want to forecast for very short-term horizon (e.g. 1 hour ahead), the NARX-RandomForest model is recommended. These two models are shown to perform consistently well on load time-series with different characteristics, without any human intervention during the training process.

%\begin{figure}
%\centering
%\label{fig:Rank}
%\includegraphics[width=0.5\textwidth,trim={0.2cm 0.3cm 0.2cm 0.2cm},clip]{Images/Rank_vs_Acc.pdf}
%\caption{Comparison between median MAPE over all 24 horizons of each method on the whole dataset}
%\end{figure}

%\begin{figure*}
%	\centering
%	\subfloat[Hvaler - Norway\label{fig:HvalerTempCons}]{
%		\includegraphics[width=.47\linewidth, trim={0cm 1.5cm 0cm 2.8cm},clip]{Images/HvalerTempCons}
%		}
%   \subfloat[GEFCom2012 - US\label{fig:Z02TempCons}]{
%		\includegraphics[width=.47\linewidth, trim={0cm 1.5cm 0cm 2.8cm},clip]{Images/Z02TempCons}
%	}
%	\caption{Relationship between Temperature ($^{\circ}\mathrm{C}$) and Consumption (KWh) in Hvaler (top) and GEFCom2012 dataset (bottom)}
%	\label{fig:TempCons}
%\end{figure*}

%%%%%%%%%%%%%%%%%%%%%%%%%%%%%%%%%%%%%%%%%%%%%%%
%%%%%%%%%%%%%%%%% CONCLUSIONS %%%%%%%%%%%%%%%%% 
%%%%%%%%%%%%%%%%%%%%%%%%%%%%%%%%%%%%%%%%%%%%%%%

\section{Conclusions and Future Work}
\label{sec:Conclusion}
In this paper, we were looking for solutions for local one-day-ahead load forecasting problem, which needs to be able to model thousands of load time-series automatically without human intervention. One baseline and five models have been proposed, including avgARIMA, orgiDSHW, modDSHW, NARX-RF, TBATS, and SemiPar. These models were tested on $40$ different load time-series, collected from US and Norway at different aggregation levels with different characteristics. The experiment results show that the SemiPar has superior performance on high-aggregation load, at the cost of a long historical data requirement. On the other hand, NARX-RF performs consistently well in many cases, at the expense of long training time. On low-aggregation load time-series, our modified version of the DSHW model works surprisingly well with only $3$ months of training data and without using temperature information. If the historical data is limited, which is the case when a new smart grid is installed, the modDSHW model is highly recommended. The experiment also suggests that at low aggregation level, long term underlying processes (e.g., trend or intra-year cycle) and temperature information do not contribute much to the forecasting accuracy. Apparently, one can develop a better and more general model for the task by automatically combine or select among those methods proposed in this paper. However, one must acknowledge the fact that this would complicate the deployment and maintenance process, where thousands of models are involved.

\bibliography{STLF_DR}
\bibliographystyle{plain}
\end{document}